\documentclass[a4paper,11pt]{article}
\usepackage[cp1251]{inputenc}
\usepackage[english]{babel}
\usepackage{amsmath,amsthm,amsfonts,amssymb}
\usepackage{graphicx}
\usepackage[left=2.5cm, right=1.0cm, top=2.0cm, bottom=2.0cm]{geometry}
\usepackage{cite}

\begin{document}

\begin{center}
\bf{Inversion formulas and their finite-dimensional analogs for\\ multidimensional Volterra equations of the first kind}

\vspace{5mm}

\it{Svetlana Solodusha, Ekaterina Antipina\\
Melentiev Energy Systems Institute SB RAS, Irkutsk, Russia}
\end{center}

\begin{abstract}
The paper focuses on  solving one class of Volterra equations of the first kind, which is characterized by the variability of all integration limits. These equations were introduced in connection with the problem of identifying nonsymmetric kernels for constructing integral models of nonlinear dynamical systems of the  "input-output"\, type in the form of Volterra polynomials. The case when the input perturbation of the system is a vector function of time is considered. To solve the identification problem, previously introduced test signals of duration $h$ (mesh step) are used in the form of linear combinations of Heaviside functions with deviating arguments. The paper  demonstrates a method for obtaining the desired solution, developing a method of steps for a one-dimensional case. The matching conditions providing the desired smoothness of the solution are established. The mesh analogs of the studied integral equations based on the formulas of middle rectangles are considered.
\end{abstract}

\section{Introduction}

Volterra integral equations are used in various applied problems. A detailed review of such applications is presented in the monograph  \cite{solodusha1}. The class of integral equations considered in this paper arises when modeling the response of a nonlinear dynamical system $y(t)$ to an input signal $x(t)$ in the form of a Volterra polynomial (a segment of an integro-power series) \cite{solodusha2}. Integral models based on Volterra polynomials attract the attention of many researchers and have an extensive field of applications (a review of the current state of research is given in  \cite{solodusha3}).  
 Of greatest interest in terms of applications is the case where the input signal $x(t)=(x_1(t),...,x_p(t))^T$ is a vector function of time:
\begin{equation} \label{Solodusha1}
y(t)=\sum\limits_{m=1}^N \sum\limits_{1\leq i_1 \leq ... \leq i_m \leq p} f_{i_1...i_m} (t),\;\; t \in [0,T],
\end{equation}
\begin{equation} \label{Solodusha2}
f_{i_1...i_m} (t) =\int\limits_0^t ... \int\limits_0^t \ K_{i_1...i_m}(s_1,..., s_m)\prod_{j=1}^m x_j(t-s_j)ds_j. 
\end{equation}
In (\ref{Solodusha1}) $y(t)$  is a scalar function of time,  $y(0)=0,$ $y'(t)\in C_{[0,T]}$. The functions  $K_{i_1...i_m}$  in  (\ref{Solodusha2}) are called Volterra kernels and they are symmetric only in those variables that correspond to the coinciding indices  $i_1...i_m$.
 The key problem in constructing a model of nonlinear dynamical system of input-output type in  form  (\ref{Solodusha1}), (\ref{Solodusha2})  lies in the identification of Volterra kernels. The absence of the symmetry property of the functions  $K_{i_1...i_m}$  complicates  the problem of constructing  (\ref{Solodusha1}).
 This paper develops the approach  \cite{solodusha4} which is based on setting    $(m-1)$-parametric families of test signals in the form of combinations of Heaviside functions with deviating arguments. 

Confine ourselves to  $x(t)=(x_1(t),x_2(t))^T$ and consider the case which is most widely used  in practice, when in  (\ref{Solodusha1}) $N=2$. Suppose further that the problem of decomposing the response   (\ref{Solodusha1}) into components  (\ref{Solodusha2})  is somehow solved and consider the problem of identifying an nonsymmetric kernel  $K_{12}$.

\section{The problem statement}

To solve the problem of identification of nonsymmetrical kernel  $K_{12}$ in
 \begin{equation} \label{Solodusha3}
f_{12}(t)=\int\limits_0^t \int\limits_0^t K_{12}(s_1,s_2)x_1(t-s_1)x_2(t-s_2)ds_1ds_2, \; \;    t \in [ 0,T ],
\end{equation}
introduce two series of  test signals  \cite{solodusha5}:
\begin{equation} \label{Solodusha4}
  \left\{\begin{array}{l}
    x_1(t)=  {e}(t)-{e}(t-h) ,\\
  x_{2_\upsilon}(t)=  {e}(t- \upsilon)-{e}(t- \upsilon -h)  ,\;t\in[h,T],\; \upsilon \leq t-h,
\end{array}\right.
\end{equation}
\begin{equation} \label{Solodusha5}
 \left\{\begin{array}{l}
 x_{1_\upsilon}(t)=   e(t- \upsilon)-{e}(t- \upsilon -h)  ,\\
  x_2(t)=     {e}(t)-{e}(t-h) ,\;t\in[h,T],\; \upsilon \leq t-h,
 \end{array}\right.
\end{equation}
where $h>0$  is  a sampling interval of the output signal, $T=Nh$, $N=$ const.  Figures 1  and  2 illustrate test signals  (\ref{Solodusha4}), (\ref{Solodusha5}), respectively.

\begin{figure}[h]
\begin{picture}(100,100)        
\put(290,-6){\vector(0,1){80}}          
\put(314,80){\makebox(0,0){$x_{2_\upsilon}(t)$}}
\put(278,6){\vector(1,0){140}}   
\put(280,40){\makebox(0,0){\small $1$}}
\put(424,-3){\makebox(0,0){$t$}}
\put(329,-4){\makebox(0,0){$\upsilon$}}
\put(365,-3){\makebox(0,0){$\upsilon+h$}}
\put(334,6){\circle*{2}}
\put(284,0){\makebox(0,0){{\small $0$}}}
\put(399,6){\circle*{2}}
\put(359,6){\circle*{2}}
\put(399,-3){\makebox(0,0){${\sc T}$}}
\put(334,35){\line(0,-4){3}}    
\put(334,25){\line(0,-4){3}}    
\put(334,15){\line(0,-4){3}}    
\put(359,35){\line(0,-4){3}}    
\put(359,25){\line(0,-4){3}}    
\put(359,15){\line(0,-4){3}}    
\thicklines
\put(334,6){\line(-4,0){44}}
\put(334,40){\line(4,0){25}}    
\put(359,6){\line(4,0){40}}    
\thinlines
\put(36,-6){\vector(0,1){80}}          
\put(26,40){\makebox(0,0){\small $1$}}
\put(55,80){\makebox(0,0){$x_{1}(t)$}}
\put(24,6){\vector(1,0){140}}   
\put(175,0){\makebox(0,0){$t$}}
\put(30,0){\makebox(0,0){{\small $0$}}}
\put(140,6){\circle*{2}}
\put(61,6){\circle*{2}}
\put(140,-3){\makebox(0,0){${\sc T}$}}
\put(61,-3){\makebox(0,0){${\sc h}$}}
\put(61,35){\line(0,-4){3}}    
\put(61,25){\line(0,-4){3}}    
\put(61,15){\line(0,-4){3}}    
\thicklines
\put(140,6){\line(-4,0){79}} 
\put(61,40){\line(-4,0){25}}   
\thinlines
\end{picture}
\caption{Test signals   $x_{1}(t)$ and $x_{2_\upsilon}(t)$.}
\end{figure}
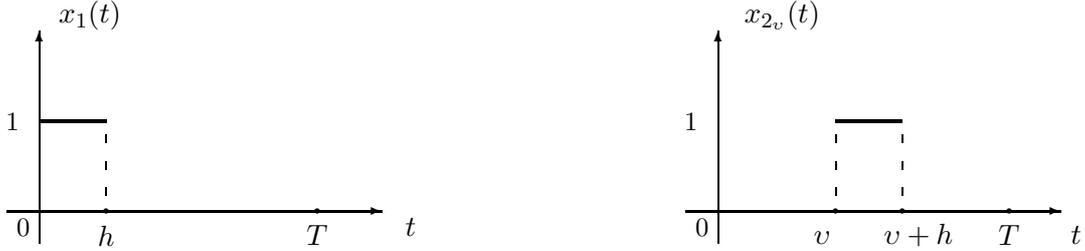

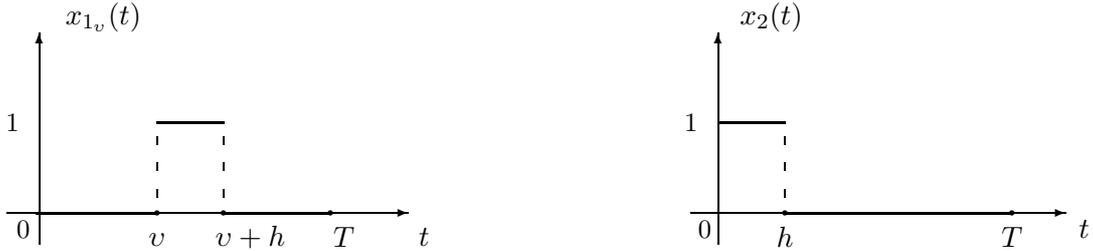
\begin{figure}[h]
\begin{picture}(100,100)                        
\put(36,-6){\vector(0,1){80}}          
\put(60,80){\makebox(0,0){$x_{1_\upsilon}(t)$}}
\put(24,6){\vector(1,0){150}}   
\put(26,40){\makebox(0,0){\small $1$}}
\put(180,-3){\makebox(0,0){$t$}}
\put(80,-4){\makebox(0,0){$\upsilon$}}
\put(115,-3){\makebox(0,0){$\upsilon+h$}}
\put(80,6){\circle*{2}}
\put(105,6){\circle*{2}}
\put(145,6){\circle*{2}}
\put(150,-3){\makebox(0,0){${\sc T}$}}
\put(80,35){\line(0,-4){3}}    
\put(80,25){\line(0,-4){3}}    
\put(80,15){\line(0,-4){3}}    
\put(105,35){\line(0,-4){3}}    
\put(105,25){\line(0,-4){3}}    
\put(105,15){\line(0,-4){3}}    
\thicklines
\put(80,40){\line(4,0){25}}   
\put(80,6){\line(-4,0){45}}  
\put(105,6){\line(4,0){40}} 
\thinlines
\put(30,0){\makebox(0,0){{\small $0$}}}
\put(290,-6){\vector(0,1){80}}          
\put(400,-3){\makebox(0,0){${\sc T}$}}
\put(315,-3){\makebox(0,0){${\sc h}$}}
\put(400,6){\circle*{2}}
\put(315,6){\circle*{2}}
\put(280,40){\makebox(0,0){\small $1$}}
\put(310,80){\makebox(0,0){$x_{2}(t)$}}
\put(280,6){\vector(1,0){140}}   
\put(427,0){\makebox(0,0){$t$}}
\put(315,35){\line(0,-4){3}}    
\put(315,25){\line(0,-4){3}}    
\put(315,15){\line(0,-4){3}}    
\put(285,0){\makebox(0,0){{\small $0$}}}
\thicklines
\put(315,40){\line(-4,0){25}}  
\put(315,6){\line(4,0){85}}  
\thinlines
\end{picture}
\caption{Test signals   $x_{1_\upsilon}(t)$  and  $x_{2}(t)$.}
\end{figure}
Substituting  (\ref{Solodusha4}), (\ref{Solodusha5})  into   (\ref{Solodusha3}) we obtain a paired Volterra equation of the first kind
\begin{equation} \label{Solodusha6}
\int\limits_{t-h}^t ds_1\int\limits^{t-\upsilon}_{t-\upsilon -h} K_{12}(s_1,s_2)ds_2=\stackrel{(1)}f(t,\upsilon),
\end{equation}
\begin{equation} \label{Solodusha7}
\int\limits^{t-\upsilon}_{t-\upsilon -h}ds_1 \int\limits_{t-h}^t K_{12}(s_1,s_2)ds_2=\stackrel{(2)}f(t,\upsilon),
\end{equation}
 where  $t\in[h,T]$,  $\stackrel{(1)}f(t,\upsilon)$ and $\stackrel{(2)}f(t,\upsilon)$  denote  responses to  (\ref{Solodusha4})   and (\ref{Solodusha5}),  respectively. It is seen from  (\ref{Solodusha6}), (\ref{Solodusha7}) that for   $\upsilon=0$,  $t\in[h,T]$  $\stackrel{(1)}f(t,0)=\stackrel{(2)}f(t,0)=f(t,0)$ holds true. 

In this case, the definition of a solution of  (\ref{Solodusha6}), (\ref{Solodusha7}) needs to be refined.   By virtue of   $t-h\geq 0$ within the low limits of integration, the domain of the sought function   $\bar{K}_{12}$ is  segment   $[0,T]$, including  segment  $[0,h]$.  Therefore, paired equation  (\ref{Solodusha6}), (\ref{Solodusha7})  makes sense only in the case where   $\bar{K}_{12}$  on  $[0,h]$  is  known.
  By  analogy  with   \cite{solodusha4},  for one-dimensional integral equations with  "prehistory", require that   
\begin{equation} \label{Solodusha8}
\bar{K}_{12}(s_1,s)={K}_{12}^{(0)}(s_1,s),\;\bar{K}_{12}(s,s_2)={K}_{12}^{(0)}(s,s_2),\;s_1,\,s_2\in[0,h),\;s\in[0,T].
\end{equation}
Consider further the procedure of obtaining the sought solution to  (\ref{Solodusha6}), (\ref{Solodusha7}).  The procedure develops the method of steps for the one-dimensional case  \cite{solodusha4}. 
 
 \section{The method of steps}
 
We  introduce   
$$\Delta_{k}=\{t,\upsilon: \, \upsilon+h\leq t,\, kh\leq t < (k+1)h\},\;     \Delta_{N+1}=\{t,\upsilon:\, \upsilon+h \leq t,\,Nh\leq t \leq T\},$$ 
$$\Delta_{0}=\{t,\upsilon: \, D_1\cup D_2,\, \upsilon\geq 0\},\; D_1=\{t,\upsilon: \, \upsilon\leq t,\,0\leq t < h,\, h>0\},$$
$$ D_2=\{t,\upsilon: \, t-h< \upsilon \leq t,\,h\leq t \leq T,\, h>0\},\;N=\frac{T}{h}, \; k=\overline{1,N},$$  such  that  
$\Delta_{0}$ coincides  with the  prehistory  and   $$\bigcup\limits_{k=1}^{N+1}\Delta_k=\{t,\upsilon: \, \upsilon+h\leq t,\,h\leq t\leq T,\,\upsilon \geq 0,\,h>0\}.$$

Let  $N(t,\upsilon)$  be  a  point  of  a  plain  with  Cartesian  coordinates.  We  will  show  that  the  condition  for  solvability  of   (\ref{Solodusha6}), (\ref{Solodusha7})  at the initial point  $N(h,0)\in\Delta_1$  is  met: 
$$ \int\limits^{t}_{t-h}ds_1\int\limits^{t}_{t-h}K_{12}(s_1,s_2)ds_2=\int\limits^{h}_{t-h}ds_1\int\limits^{h}_{t-h}K_{12}(s_1,s_2)ds_2+\int\limits^{h}_{t-h}ds_1\int\limits^{t}_{h}K_{12}(s_1,s_2)ds_2+$$
$$+\int\limits^{t}_{h}ds_1\int\limits^{h}_{t-h}K_{12}(s_1,s_2)ds_2+\int\limits^{t}_{h}ds_1\int\limits^{t}_{h}K_{12}(s_1,s_2)ds_2=\int\limits^{h}_{t-h}ds_1\int\limits^{h}_{t-h}K_{12}^{(0)}(s_1,s_2)ds_2+$$ $$+\int\limits^{h}_{t-h}ds_1\int\limits^{t}_{h}K_{12}^{(0)}(s_1,s_2)ds_2+\int\limits^{t}_{h}ds_1\int\limits^{h}_{t-h}K_{12}^{(0)}(s_1,s_2)ds_2+\int\limits^{t}_{h}ds_1\int\limits^{t}_{h}K_{12}(s_1,s_2)ds_2
$$
and   by   (\ref{Solodusha8})
$$ \int\limits^{t}_{h}ds_1\int\limits^{t}_{h}K_{12}(s_1,s_2)ds_2=f(t,0)-\int\limits^{h}_{t-h}ds_1\int\limits^{h}_{t-h}K_{12}^{(0)}(s_1,s_2)ds_2-$$
$$-\int\limits^{h}_{t-h}ds_1\int\limits^{t}_{h}K_{12}^{(0)}(s_1,s_2)ds_2-\int\limits^{t}_{h}ds_1\int\limits^{h}_{t-h}K_{12}^{(0)}(s_1,s_2)ds_2\equiv f_1(t,0). $$

Note,  that 
$$f_1(h,0)=f(h,0)- \int\limits^{h}_{0}ds_1\int\limits^{h}_{0}K_{12}^{(0)}(s_1,s_2)ds_2\equiv 0.$$
Assuming   that  $\stackrel{(1)}f(t,\upsilon),\,\stackrel{(2)}f(t,\upsilon) \in C^{(2)}_{\Delta_1}$, we  solve    (\ref{Solodusha6}), (\ref{Solodusha7}) by differentiation  with  respect  to   $t$  and  $\upsilon$,
so   that   for   $N(t,\upsilon)\in\Delta_1$
\begin{equation} \label{Solodusha9}
\bar{K}_{12}(M)={\cal D}_2\stackrel{(1)}{f}(t,\upsilon)
+{K}_{12}^{(0)}(t,t-\upsilon-h)
+{K}_{12}^{(0)}(t-h,t-\upsilon)
-{K}_{12}^{(0)}(t-h,t-\upsilon-h),\end{equation}
\begin{equation} \label{Solodusha10}
\bar{K}_{12}({\bar M})={\cal D}_2\stackrel{(2)}{f}(t,\upsilon)
+{K}_{12}^{(0)}(t-\upsilon-h,t)
+{K}_{12}^{(0)}(t-\upsilon,t-h)
-{K}_{12}^{(0)}(t-\upsilon-h,t-h),\end{equation}
$${\cal D}_2\stackrel{(1)}{f}(t,\upsilon)=-\left(\stackrel{(1)}{f''_{t\upsilon}}+\stackrel{(1)}{f''_{\upsilon^2}}\right),\;{\cal D}_2\stackrel{(2)}{f}(t,\upsilon)=-\left(\stackrel{(2)}{f''_{t\upsilon}}+\stackrel{(2)}{f''_{\upsilon^2}}\right),$$
where  $M(p,q)\in \Omega_k(N(t,\upsilon))$, $k=\overline{1,N},$   is  a  point  of  the  plain  with  Cartesian  coordinates   $(p,q)$,  $t-\upsilon\leq q \leq p\leq t$, $\bar{M}$  is  a  point  from   $\Omega_k(N(t,\upsilon))$, symmetrical  to  $M$  with  respect  to  diagonal  $p=q$,  so   that  $t-\upsilon\leq p \leq q \leq t$. 

We  denote  this  solution  by  $K_{12}^{(1)}$  and  rewrite   (\ref{Solodusha9}), (\ref{Solodusha10})  in  the  following  form:
\begin{equation} \label{Solodusha010}
{K}_{12}^{(i)}(M)={\cal D}_2\stackrel{(1)}{f}\big|_{N(t,\upsilon)}
+{K}_{12}^{(i-1)}(t,t-\upsilon-h)
+{K}_{12}^{(i-1)}(t-h,t-\upsilon)
-{K}_{12}^{(i-1)}(t-h,t-\upsilon-h),
\end{equation}
\begin{equation} \label{Solodusha011}
{K}_{12}^{(i)}(\bar{M})={\cal D}_2\stackrel{(2)}{f}\big|_{N(t,\upsilon)}
+{K}_{12}^{(i-1)}(t-\upsilon-h,t)
+{K}_{12}^{(i-1)}(t-\upsilon,t-h)
-{K}_{12}^{(i-1)}(t-\upsilon-h,t-h),
\end{equation}
where   $i=1$.

The condition of simultaneous continuity of the initial function  ${K}_{12}^{(0)}$ and  the  desired solution  ${K}_{12}^{(1)}$  at  points  $M,$ $\bar{M}\in \Omega_1(N(t,\upsilon))$  for  $h \leq t < 2h$, $ \upsilon= t-h$  follows  from  (\ref{Solodusha010}),  (\ref{Solodusha011}):
$$
{\cal D}_2\stackrel{(1)}{f}\big|_{N(t,\upsilon)}={K}_{12}^{(i)}(t,t-\upsilon)-
{K}_{12}^{(i-1)}(t,t-\upsilon-h) -{K}_{12}^{(i-1)}(t-h,t-\upsilon)
+{K}_{12}^{(i-1)}(t-h,t-\upsilon-h)=
$$
\begin{equation} \label{Solodusha0010}={K}_{12}^{(i-1)}(t,t-\upsilon)-
{K}_{12}^{(i-1)}(t,t-\upsilon-h) -{K}_{12}^{(i-1)}(t-h,t-\upsilon)
+{K}_{12}^{(i-1)}(t-h,t-\upsilon-h),\end{equation}
$$
{\cal D}_2\stackrel{(2)}{f}\big|_{N(t,\upsilon)}={K}_{12}^{(i)}({t-\upsilon,t})-
{K}_{12}^{(i-1)}(t-\upsilon-h,t)-{K}_{12}^{(i-1)}(t-\upsilon,t-h)
+{K}_{12}^{(i-1)}(t-\upsilon-h,t-h)=
$$
\begin{equation} \label{Solodusha0011}
={K}_{12}^{(i-1)}({t-\upsilon,t})-
{K}_{12}^{(i-1)}(t-\upsilon-h,t)-{K}_{12}^{(i-1)}(t-\upsilon,t-h)
+{K}_{12}^{(i-1)}(t-\upsilon-h,t-h).
\end{equation}
In particular,  from   (\ref{Solodusha0010}), (\ref{Solodusha0011}) 
  for  $M(h,h)={\bar M}(h,h)$ ($t=h,$ $\upsilon=0$),  we  have 
$$ {\cal D}_2{f}\big|_{N(h,0)}=
{K}_{12}^{(0)}(h,h)-{K}_{12}^{(0)}(h,0)-{K}_{12}^{(0)}(0,h)+{K}_{12}^{(0)}(0,0),$$
$${\cal D}_2{f}\big|_{N(h,0)}={\cal D}_2\stackrel{(1)}{f}\big|_{N(h,0)}={\cal D}_2\stackrel{(2)}{f}\big|_{N(h,0)}.$$
  If   ${K}_{12}^{(0)}$   is  continuous  on  $\Omega_0,$   then (\ref{Solodusha010}), (\ref{Solodusha011})  implies  that  ${K}_{12}^{(1)}$  is  continuous  on  $\Omega_1$. 

 Now let  $N(t,0)\in \Delta_2$.  Then,  using  the  same  procedure,  we  have 
$$\int\limits_{t-h}^{t} ds_1 \int\limits_{t-h}^{t} K_{12}(s_1,s_2)ds_2= \int\limits_{t-h}^{2h} ds_1 \int\limits_{t-h}^{2h} K_{12}(s_1,s_2)ds_2+
\int\limits_{t-h}^{2h} ds_1 \int\limits_{2h}^{t} K_{12}(s_1,s_2)ds_2+
$$
$$
+\int\limits_{2h}^{t} ds_1 \int\limits_{t-h}^{2h} K_{12}(s_1,s_2)ds_2+
\int\limits_{2h}^{t} ds_1 \int\limits_{2h}^{t} K_{12}(s_1,s_2)ds_2=\int\limits_{t-h}^{2h} ds_1 \int\limits_{t-h}^{2h} K_{12}^{(1)}(s_1,s_2)ds_2+
$$
$$
+
\int\limits_{t-h}^{2h} ds_1 \int\limits_{2h}^{t} K_{12}^{(1)}(s_1,s_2)ds_2+
\int\limits_{2h}^{t} ds_1 \int\limits_{t-h}^{2h} K_{12}^{(1)}(s_1,s_2)ds_2+
\int\limits_{2h}^{t} ds_1 \int\limits_{2h}^{t} K_{12}(s_1,s_2)ds_2,
$$
so that (\ref{Solodusha6}),  (\ref{Solodusha7})  implies
\begin{equation} \label{Solodusha12}
\int\limits_{2h}^{t} ds_1 \int\limits_{2h}^{t} K_{12}(s_1,s_2)ds_2=f(t,0)-\int\limits_{t-h}^{2h} ds_1 \int\limits_{t-h}^{2h} K_{12}^{(1)}(s_1,s_2)ds_2-\end{equation}
$$
-\int\limits_{t-h}^{2h} ds_1 \int\limits_{2h}^{t} K_{12}^{(1)}(s_1,s_2)ds_2-\int\limits_{2h}^{t} ds_1 \int\limits_{t-h}^{2h} K_{12}^{(1)}(s_1,s_2)ds_2\equiv f_2(t,0).
$$

The  solvability  condition  for  (\ref{Solodusha12}) $f_2(2h,0)=0$ is obviously satisfied, since
$$
f_2(2h,0)=f(2h,0)-\int\limits_{h}^{2h} ds_1 \int\limits_{h}^{2h} K_{12}^{(1)}(s_1,s_2)ds_2=$$
$$=\int\limits_{h}^{2h} ds_1 \int\limits_{h}^{2h} K_{12}^{(1)}(s_1,s_2)ds_2-\int\limits_{h}^{2h} ds_1 \int\limits_{h}^{2h} K_{12}^{(1)}(s_1,s_2)ds_2\equiv 0.
$$
Therefore,  under the  assumption  that   $f_2(t,\upsilon)\in C_{\Delta_2}^{(2)}$   we  have   (\ref{Solodusha010}),  (\ref{Solodusha011}),
where  $N(t,\upsilon)\in {\Delta_2}$, $i=2$.

Make  sure  that  conditions  (\ref{Solodusha0010}),  (\ref{Solodusha0011})  for  $i=2$ provide continuity,  firstly, of functions   $K_{12}^{(2)}$  and  $K_{12}^{(0)}$  at  $2h \leq t<3h$, $\upsilon=t-h$, and, secondly, of functions  $K_{12}^{(1)}$,  $K_{12}^{(2)}$  at  $t=2h$, $0 \leq \upsilon< h$. Indeed, in the first case it follows from  (\ref{Solodusha0010})  that   
$$ 
{\cal D}_2\stackrel{(1)}{f}\big|_{N(t,\upsilon)}={K}_{12}^{(2)}(t,h)-
{K}_{12}^{(1)}(t,0) -{K}_{12}^{(1)}(t-h,h)
+{K}_{12}^{(1)}(t-h,0)={K}_{12}^{(2)}(t,h)-
$$
$$
-{K}_{12}^{(0)}(t,0) -{K}_{12}^{(0)}(t-h,h)
+{K}_{12}^{(0)}(t-h,0)={K}_{12}^{(0)}(t,h)
-{K}_{12}^{(0)}(t,0) -{K}_{12}^{(0)}(t-h,h)
+{K}_{12}^{(0)}(t-h,0),
$$
similarly, from  (\ref{Solodusha0011}) we have 
$${\cal D}_2\stackrel{(2)}{f}\big|_{N(t,\upsilon)}={K}_{12}^{(0)}({h,t})-
{K}_{12}^{(0)}(0,t)-{K}_{12}^{(0)}(h,t-h)
+{K}_{12}^{(0)}(0,t-h),
$$
so  that  $$\lim_{\varepsilon\rightarrow 0} {K}_{12}^{(0)}(t,h-\varepsilon) ={K}_{12}^{(2)}(t,h),\;\lim_{\varepsilon\rightarrow 0} {K}_{12}^{(0)}(h-\varepsilon,t) ={K}_{12}^{(2)}(h,t).$$
In the second case, the substitution of $i=2$, $t=2h$ in (\ref{Solodusha0010}), (\ref{Solodusha0011})  gives 
$$
{\cal D}_2\stackrel{(1)}{f}\big|_{N(2h,\upsilon)}={K}_{12}^{(1)}(2h,2h-\upsilon)
+{K}_{12}^{(1)}(2h,h-\upsilon)
+{K}_{12}^{(1)}(h,2h-\upsilon)
-{K}_{12}^{(1)}(h,h-\upsilon),
$$
$$
{\cal D}_2\stackrel{(2)}{f}\big|_{N(2h,\upsilon)}={K}_{12}^{(1)}(2h-\upsilon,2h)
+{K}_{12}^{(1)}(h-\upsilon,2h)
+{K}_{12}^{(1)}(2h-\upsilon,h)
-{K}_{12}^{(1)}(h-\upsilon,h),
$$
so  that   $$\lim_{\varepsilon\rightarrow 0} {K}_{12}^{(1)}(2h-\varepsilon,2h-\upsilon-\varepsilon) ={K}_{12}^{(2)}(2h,2h-\upsilon),\;\lim_{\varepsilon\rightarrow 0} {K}_{12}^{(1)}(2h-\upsilon-\varepsilon,2h-\varepsilon) ={K}_{12}^{(2)}(2h-\upsilon,2h).$$

Extending  by  \cite{solodusha4}  this  process  to    $\Delta_k,$ $k=\overline{1,N+1}$, we  find  a  solution  to  (\ref{Solodusha6}), (\ref{Solodusha7}):
$$\bar{K}_{12}(t,t-\upsilon)=\sum_{i=1}^{N+1}{\cal D}_2\stackrel{(1)}{f}\big|_{N\in \Delta_i}+ \sum_{i=1}^{N+1} \left(
{K}_{12}^{(i-1)}(t,t-\upsilon-h)+\right.$$
$$\left.
+{K}_{12}^{(i-1)}(t-h,t-\upsilon)
-{K}_{12}^{(i-1)}(t-h,t-\upsilon-h)\right),
$$
$$
\bar{K}_{12}(t-\upsilon,t)=\sum_{i=1}^{N+1}{\cal D}_2\stackrel{(2)}{f}\big|_{N\in \Delta_i}
+ \sum_{i=1}^{N+1} \left({K}_{12}^{(i-1)}(t-\upsilon-h,t)+\right.$$
$$\left.
+{K}_{12}^{(i-1)}(t-\upsilon,t-h)
-{K}_{12}^{(i-1)}(t-\upsilon-h,t-h)\right).$$

\section{Numerically  solving  (\ref{Solodusha6}), (\ref{Solodusha7})}

Let us solve numerically   (\ref{Solodusha6}), (\ref{Solodusha7}). For simplicity,  denote  $K_{12}(t_{i-\frac{1}{2}},\,t_{i-j-\frac{1}{2}})$  by $K_{i-\frac{1}{2},\,i-j-\frac{1}{2}}$. 
Introduce for $h\leq t \leq T$ a uniform mesh by coordinating integer nodes with the points of discontinuities on (\ref{Solodusha4}),  (\ref{Solodusha5})
$$
t_i=ih,\; t_{i-\frac{1}{2}}=\left( i-\frac{1}{2}\right)h, \; \upsilon=jh, \; j=\overline{0,i-2}, \; i=\overline{2,N}, \; T=Nh.
$$
Applying the quadrature formula of the middle rectangles for the approximation of the integrals in the left-hand side of  (\ref{Solodusha6}), (\ref{Solodusha7}) we  find  that the discrete analog of  (\ref{Solodusha6}), (\ref{Solodusha7}) has  the  form
\begin{equation} \label{Solodusha16}
h^2 \sum \limits_{l=i}^i \sum\limits_{m=i-j}^{i-j} K^h_{l-\frac{1}{2},m-\frac{1}{2}}=\stackrel{(1)}{f^h}_{i,j},\;    j=\overline{0,i-2},
\end{equation}
\begin{equation} \label{Solodusha17}
h^2 \sum\limits_{l=i-j}^{i-j} \sum \limits_{m=i}^i K^h_{l-\frac{1}{2},m-\frac{1}{2}}=\stackrel{(2)}{f^h}_{i,j},\;     j=\overline{1,i-2}
\end{equation}
 ((\ref{Solodusha17}) factors in the equality   $\stackrel{(1)}{f^h_{i,0}}=\stackrel{(2)}{f^h_{i,0}}$  for  $i=\overline{1,N}$). 
 Inverse  formulas  of  a  paired  system  of  linear  algebraic  equations (\ref{Solodusha16}),  (\ref{Solodusha17})  have  the  form  
$$
K_{i-\frac{1}{2},i-j-\frac{1}{2}}=\frac{\stackrel{(1)}{f^h}(i,j)}{h^2},\;   j=\overline{0,i-2},\;\;
K_{i-j-\frac{1}{2},i-\frac{1}{2}}=\frac{\stackrel{(2)}{f^h}(i,j)}{h^2},\;   j=\overline{1,i-2}.
$$
We will demonstrate numerical calculations using  (\ref{Solodusha16}), (\ref{Solodusha17}) on  an  example.

Let  an  exact  solution  to  (\ref{Solodusha6}), (\ref{Solodusha7})  be  the  function  
$$
K(s_1,s_2)=as_1^2-bs_2,\; a,b=\mbox{const}.
$$
Then  the  right-hand  sides  of   (\ref{Solodusha6}), (\ref{Solodusha7})  have  the  form  
$$
\stackrel{(1)}f(t,\upsilon)= \frac{ah^2}{3}  \left( 3t(t-h)+h^2  \right)  +\frac{bh^2}{2}(h-2(t-\upsilon)),
$$
$$
\stackrel{(2)}f(t,\upsilon)=\frac{ah^2}{3}  \left(3 (t - \upsilon)(t- \upsilon-h) +h^2 \right)  +\frac{bh^2}{2}(h-2t).
$$
Determine the difference approximation  $K^h$  using  (\ref{Solodusha16}), (\ref{Solodusha17})  and  find  errors  
$$
\varepsilon_1=\max_{i,\,j} \vert  K^h_{i-\frac{1}{2},i-j-\frac{1}{2}}-K(t_{i-\frac{1}{2}},t_{i-j-\frac{1}{2}})   \vert,\; \varepsilon_2=\max_{i,\,j} \vert  K^h_{i-j-\frac{1}{2},i-\frac{1}{2}}-K(t_{i-j-\frac{1}{2}},t_{i-\frac{1}{2}})   \vert.
$$
Table 1 presents the results of the calculations of   $\varepsilon=\max\{\varepsilon_1,\varepsilon_2 \}$  for  $a=4,$ $b=-1$. 
\begin{center}
\begin{table}[h]
\caption{The quadrature method of middle rectangles.}
\centering
\begin{tabular}{|c|c|}
\hline
$h$ & $\varepsilon$\\
\hline
0{.}2500  & 0{.}00520  \\
\hline 
0{.}1250  & 0{.}00130  \\
\hline
0{.}0625  & 0{.}00033 \\
\hline
\end{tabular}
\end{table}
\end{center}

As can be seen from the Table, the numerical method has the second order of convergence, i.e.  with a decrease in the mesh step by half,  $\varepsilon$ decreases by a factor of 4.

\section{Conclusion}

The paper considers solving the paired two-dimensional Volterra integral equation of the first kind arising in the problem of identification of nonsymmetric Volterra kernels. The method of obtaining the desired solution develops the method of steps for the one-dimensional case. The coordination conditions that ensure the continuity of the solution are indicated. A mesh analog of the solution obtained on the basis of cubes of middle rectangles is given. 

\subsection*{Acknowledgments}
The research was carried out under State Assignment III.17.3.1  of  the Fundamental Research of Siberian Branch of the Russian Academy of Sciences, reg. No. AAAA-A17-117030310442-8.


\section*{References}
\begin{thebibliography}{9}
\bibitem{solodusha1} Brunner H 2017 {\it Volterra integral equations: an introduction to theory and applications} (Cambridge:  Cambridge University Press)
\bibitem{solodusha2} Volterra V 1959 {\it Theory of Functionals and of Integral and Integro-Differential Equations} (New York: Dover  Publications) 
\bibitem{solodusha3} Cheng C M, Peng Z K, Zhang W M and Meng G 2017 {Volterra-series-based nonlinear system modeling and its engineering applications: A state-of-the-art review} {\it Mechanical Systems and Signal Processing}  {\bf{87}}  340--64
\bibitem{solodusha4} Apartsyn A S 2003 {\it Nonclassical Linear Volterra Equations of the First Kind} (Utrecht: VSP)
\bibitem{solodusha5} Solodusha S V 2019 {\it Methods for constructing integral models of dynamic systems: algorithms and applications in power engineering} (Irkutsk: ESI SB RAS)
\end{thebibliography}
\end{document}